\documentclass[secthm,seceqn,amsthm,ussrhead,12pt]{amsart}
\usepackage{amsmath,latexsym}
\usepackage[english]{babel}
\usepackage[psamsfonts]{amssymb}
\usepackage{times}
\usepackage{cite}
\usepackage{pdflscape} 
\usepackage{ulem}
\usepackage[mathcal]{euscript}
\usepackage{tikz}
\usepackage{hyperref}
\usepackage{colortbl}
\usepackage{cancel}
\usepackage{color}

\usetikzlibrary{arrows}

\newtheorem{Th}{Theorem}
\newtheorem{Lem}[Th]{Lemma}

\newtheorem{corollary}[Th]{Corollary}

\newenvironment{Proof}[1][Proof.]{\begin{trivlist}
\item[\hskip \labelsep {\bfseries #1}]}{\flushright
$\Box$\end{trivlist}}

\begin{document}
	\sloppy


{\Large Degenerations of Jordan Algebras and ''Marginal'' Algebras
\footnote{The work was supported by   
FAPESP 	16/16445-0, 	18/15712-0; 
RFBR  18-31-00001;
and the President's "Program Support of Young Russian Scientists" (grant MK-2262.2019.1).}
}

\medskip

\medskip

\medskip

\medskip
\textbf{Ilya Gorshkov$^{a}$, Ivan Kaygorodov$^{b}$, Yury Popov$^{c}$}
\medskip

{\tiny

$^{a}$ Sobolev Institute of Mathematics, Novosibirsk, Russia.

$^{b}$ Universidade Federal do ABC, CMCC, Santo Andr\'{e}, Brazil.

$^{c}$ Universidade Estadual de Campinas, Campinas, Brazil.

\smallskip

    E-mail addresses:\smallskip

    Ilya Gorshkov (ilygor8@gmail.com),

    \smallskip

    Ivan Kaygorodov (kaygorodov.ivan@gmail.com),

    \smallskip
 
    Yury Popov (yuri.ppv@gmail.com).    
\smallskip

}

\

{\bf Abstract.} 
We describe all degenerations of the variety $\mathfrak{Jord}_3$ of Jordan algebras of dimension three over $\mathbb{C}.$
In particular, we describe all irreducible components in $\mathfrak{Jord}_3.$  
For every $n$ we define an $n$-dimensional rigid ''marginal'' Jordan algebra of level one.
Also, we discuss about ''marginal'' algebras in associative, alternative, left alternative, non commutative Jordan, Leibniz and anticommutative cases.

\

       \vspace{0.3cm}


\section{Introduction}


Degenerations of algebras is an interesting subject that was studied in various papers 
(see, for example, \cite{ale,ale2,aleis,maria,gkks,kpv19,kkk18,ack,BC99,S90,GRH,GRH2,BB14,ikv17,ikv18}).
In particular, there are many results concerning degenerations of algebras of low dimensions in a  variety defined by a set of identities.
One of important problems in this direction is the description of so-called rigid algebras. These algebras are of big interest, since the closures of their orbits under the action of generalized linear group form irreducible components of a variety under consideration
(with respect to the Zariski topology). For example, the rigid algebras were classified in the varieties of
low dimensional Jordan  \cite{KS07,K1} and Leibniz \cite{ikv17} algebras.
There are fewer works in which the full information about degenerations was found for some variety of algebras.
This problem was solved 
for two-dimensional Jordan algebras in \cite{jor2},
for two-dimensional Pre-Lie algebras in \cite{BB09},
for three-dimensional Jordan superalgebras in \cite{maria},
for four-dimensional Lie algebras in \cite{BC99}, 
for nilpotent five- and six-dimensional Lie algebras in \cite{S90,GRH}, 
for three-dimensional Novikov algebras in \cite{BB14}, 
for nilpotent five- and six-dimensional Malcev algebras in \cite{kpv},
for four-dimensional Zinbiel and nilpotent Leibniz algebras \cite{kpv16},
and all $2$-dimensional algebras \cite{kv16}.
Interesting notion concerning degenerations is the so-called level of an algebra, which is the maximal length of a chain of nontrivial degenerations starting with this algebra. 
The algebras of the first level  are classified in \cite{khud13}.
Another interesting notion concerning degenerations is the study of irreducible components (see, for example, \cite{khom2}).
In the paper  \cite{gorb12}, the author proved that in the variety on $n$-dimensional Lie algebras there exists an algebra with nontrivial multiplication which lies in all irreducible components.

Jordan algebras appeared as a tool for studies in quantum mechanic in
the paper of Jordan, von Neumann and Wigner \cite{jnw}.
A commutative algebra is called a {\it  Jordan  algebra}  if it satisfies the identity
$$(x^2y)x=x^2(yx).$$ 
The study of the structure theory and other properties of Jordan algebras was initiated by Albert in \cite{albert}.

The paper is organized as follows.
In the sections \ref{def} and \ref{met} we give the definitions, notation and methods that we use in the paper.  
The section \ref{deg} contains full information about degenerations of Jordan algebras of dimension $3.$
Particularly, we construct the graph of primary degenerations. 
The vertices of this graph are the isomorphism classes of algebras in the variety under consideration.
An algebra $A$ degenerates to an algebra $B$ iff there is a path from the vertex corresponding to $A$ to the vertex corresponding to $B$.
Also we describe rigid algebras and irreducible components in the varieties of algebras under consideration.
The section \ref{marginal} is dedicated to a study of ''marginal'' algebras. 
In this section we prove that for each dimension $n \geq 2$ there exists an $n$-dimensional rigid Jordan algebra of level one. Hence, it follows that there is no Jordan $n$-dimensional algebra, $n \geq 2$, with nonzero multiplication which lies in the intersection of all irreducible components in the variety of $n$-dimensional Jordan algebras.
Therefore, we show that the analouge of result of Gorbatsevich \cite{gorb12} is not true for Jordan algebras, 
associative, alternative, left alternative and Leibniz algebras
and it is true for Malcev, binary Lie, Sagle and other anticommutative algebras.

$\mathfrak{Remark.}$
Note that some particular results of this paper were obtained in \cite{me}.
In particular, some elementary non-degenerations were proved and it was proved that there is no algebra which degenerates to the ''marginal''  Jordan algebra.
In this paper we give the complete description of all degenerations and non-degenerations in a more elegant and short form and 
prove that the marginal Jordan algebra is rigid.
It is a generalization of the results obtained in \cite{me}.


\section{Definitions and notation}\label{def}

All spaces in this paper are considered over $\mathbb{C}$, and we write simply $dim$, $Hom$ and $\otimes$ instead of $dim_{\mathbb{C}}$, $Hom_{\mathbb{C}}$ and $\otimes_{\mathbb{C}}$. 
An algebra $\mathbb{A}$ is a set with a structure of a vector space and a linear map from $\mathbb{A}\otimes \mathbb{A}$ to $\mathbb{A}$.

Given an $n$-dimensional vector space $\mathbb{V}$, the set $Hom(\mathbb{V} \otimes \mathbb{V},\mathbb{V}) \cong \mathbb{V}^* \otimes \mathbb{V}^* \otimes \mathbb{V}$ 
is a vector space of dimension $n^3$. This space has a structure of the affine variety $\mathbb{C}^{n^3}.$ Indeed, let us fix a basis $e_1,\dots,e_n$ of $\mathbb{V}$. Then any $\mu\in Hom(\mathbb{V} \otimes \mathbb{V},\mathbb{V})$ is determined by $n^3$ structure constants $c_{i,j}^k\in\mathbb{C}$ such that
$\mu(e_i\otimes e_j)=\sum\limits_{k=1}^nc_{i,j}^ke_k$. A subset of $Hom(\mathbb{V} \otimes \mathbb{V},\mathbb{V})$ is {\it Zariski-closed} if it can be defined by a set of polynomial equations in the variables $c_{i,j}^k$ ($1\le i,j,k\le n$).

Let $\mathfrak{Jord}_n$ be the set of all Jordan algebras of dimension $n,$ understood as a subset of an affine variety $Hom(\mathbb{V} \otimes \mathbb{V}, \mathbb{V}).$ Then one can see that $\mathfrak{Jord}_n$ is a Zariski-closed subset of the variety $Hom(\mathbb{V} \otimes \mathbb{V},\mathbb{V})$ (the defining identities of Jordan algebras induce corresponding polynomial relations in the variables $c_{i,j}^k$). 
The general linear group ${\rm GL}(\mathbb{V})$ acts on $\mathfrak{Jord}_n$ by conjugations:
$$ (g * \mu )(x\otimes y) = g\mu(g^{-1}x\otimes g^{-1}y)$$ 
for $x,y\in V$, $\mu\in \mathfrak{Jord}_n \subset Hom(\mathbb{V} \otimes \mathbb{V},\mathbb{V})$ and $g\in {\rm GL}(\mathbb{V})$.
Thus, $\mathfrak{Jord}_n$ is decomposed into ${\rm GL}(\mathbb{V})$-orbits that correspond to the isomorphism classes of algebras. 
Let $O(\mu)$ denote the orbit of $\mu\in \mathfrak{Jord}_n$ under the action of ${\rm GL}(\mathbb{V})$ and $\overline{O(\mu)}$ denote the Zariski closure of $O(\mu)$.

Let $A$ and $B$ be two $n$-dimensional Jordan algebras and let $\mu,\lambda \in \mathfrak{Jord}_n$ represent $A$ and $B$ respectively.
We say that $A$ degenerates to $B$ and write $A\to B$ if $\lambda\in\overline{O(\mu)}$.
Note that in this case we have $\overline{O(\lambda)}\subset\overline{O(\mu)}$. 
Hence, the definition of a degeneration does not depend on the choice of $\mu$ and $\lambda$. 
If $A\not\cong B$, then the assertion $A\to B$ is called a {\it proper degeneration}. We write $A\not\to B$ if $\lambda\not\in\overline{O(\mu)}$.

Let $A$ be represented by $\mu\in \mathfrak{Jord}_n.$ Then  $A$ is  {\it rigid} in $\mathfrak{Jord}_n$ if $O(\mu)$ is an open subset of $\mathfrak{Jord}_n$. One of the ways (which has the advantage to be independent of classification) of proving the rigidity of an algebra $\mathfrak{J}$ is to show that its second cohomology group $H^2(\mathfrak{J},\mathfrak{J})$ with coefficients in itself is zero (the definition of the second cohomology group can be seen in the section \ref{mar}).

 Recall that a subset of a variety is called irreducible if it cannot be represented as a union of two non-trivial closed subsets. A maximal irreducible closed subset of a variety is called its {\it irreducible component}.
 In particular, $A$ is rigid in $\mathfrak{Jord}_n$ iff $\overline{O(\mu)}$ is an irreducible component of $\mathfrak{Jord}_n$. 
 It is well known that any affine variety can be represented as a finite union of its irreducible components in a unique way. We denote by $Rig(\mathfrak{Jord}_n)$ the set of rigid algebras in $\mathfrak{Jord}_n.$
We use the notation $S_i=\langle e_i,\dots,e_n\rangle,\ i=1,\ldots,n$.



\section{Methods}\label{met} 

In the present work we use the methods for proving degenerations and non-degenerations applied to Lie and Jordan algebras in \cite{BC99,GRH,GRH2,S90,KS07}.
First of all, if $A\to B$ and $A\not\cong B$, then $\mathfrak{Der}(A)<\mathfrak{Der}(B)$, where $\mathfrak{Der}(A)$ is the algebra of derivations of $A$. We compute the dimensions of algebras of derivations and  check the assertion $A\to B$ only for such $A$ and $B$ that $\mathfrak{Der}(A)<\mathfrak{Der}(B)$. Secondly, if $A\to C$ and $C\to B$ then $A\to B$. If there is no $C$ such that $A\to C$ and $C\to B$ are proper degenerations, then the assertion $A\to B$ is called a {\it primary degeneration}. If $\mathfrak{Der}(A)<\mathfrak{Der}(B)$ and there are no $C$ and $D$ such that $C\to A$, $B\to D$, $C\not\to D$ and one of the assertions $C\to A$ and $B\to D$ is a proper degeneration,  then the assertion $A \not\to B$ is called a {\it primary non-degeneration}. It suffices to prove only primary degenerations and non-degenerations to describe degenerations in the variety under consideration. It is easy to see that any algebra degenerates to the algebra with zero multiplication. From now on we use this fact without mentioning it.

To prove primary degenerations, we will construct families of matrices parametrized by $t$. Namely, let $A$ and $B$ be two algebras represented by the structures $\mu$ and $\lambda$ from $\mathfrak{Jord}_n$ respectively. Let $e_1,\dots, e_n$ be a basis of $\mathbb{V}$ and $c_{i,j}^k$ ($1\le i,j,k\le n$) be the structure constants of $\lambda$ in this basis. If there exist rational functions $a_{ij}(t)\in\mathbb{C}$ ($1\le i,j\le n$, $t\in\mathbb{C}^*$) such that the vectors $E_i^t=\sum\limits_{j=1}^na_{ij}(t)e_j$ ($1\le i\le n$) form a basis of $\mathbb{V}$ for any $t\in\mathbb{C}^*$, and the structure constants of $\mu$ in the basis $E_1^t,\dots, E_n^t$ are such rational functions $c_{i,j}^k(t)\in\mathbb{C}[t]$ that $c_{i,j}^k(0)=c_{i,j}^k$, then $A\to B$. In this case  $E_1^t,\dots, E_n^t$ is called a {\it parametrized basis} for $A\to B$.


Let us describe the methods for proving primary non-degenerations. The main tool for this is the following lemma.

\begin{Lem}[\cite{GRH}]\label{main}
Let $\mathcal{B}$ be a Borel subgroup of ${\rm GL}(\mathbb{V})$ and $\mathcal{R}\subset \mathfrak{Jord}_n$ be a $\mathcal{B}$-stable closed subset.
If $A \to B$ and $A$ can be represented by $\mu\in\mathcal{R}$ then there is $\lambda\in \mathcal{R}$ that represents $B$.
\end{Lem}


To prove the degeneration $A\not\to B$
we will define $\mathcal{R}$ by a set of polynomial equations and will give a basis of $\mathbb{V}$, in which the structure constants of $\mu$ give a solution to all these equations. We will omit everywhere the verification of the fact that $\mathcal{R}$ is stable under the action of the subgroup of upper triangular matrices and of the fact that $\lambda\not\in\mathcal{R}$ for any choice of a basis of $\mathbb{V}$. These verifications can be done by direct calculations.

If the number of orbits under the action of ${\rm GL}(\mathbb{V})$ on  $\mathfrak{Jord}_n$ is finite, then the graph of primary degenerations gives the whole picture. 
In particular, the description of rigid algebras and irreducible components can be easily obtained. 

The degenerations must also preserve the action of semisimple subalgebras:
\begin{Lem}[\cite{K1}]
\label{ss_action}
For a non-nilpotent Jordan algebra $\mathfrak{J}$ consider the pair $(\mathfrak{J}_{ss}, \Gamma(\mathfrak{J}))$, where $\mathfrak{J}_{ss}$ is the maximal semisimple subalgebra of $\mathfrak{J}$ and $\Gamma(\mathfrak{J})$ is the representation of $\mathfrak{J}_{ss}$ on $\mathfrak{J}$ (obtained by the restriction of the regular representation of $\mathfrak{J}$). Then for any irreducible component $\mathbb{T}$ of $\mathfrak{Jor}_n$ , $(\mathfrak{J}_{ss}, \Gamma(\mathfrak{J}))$ is constant on an open
subset of $\mathbb{T}.$ If $\mathfrak{J}\to \mathfrak{J}_1 \in \mathbb{T}$ then $(\mathfrak{J}_1)_{ss}$ is a subalgebra of $\mathfrak{J}_{ss}$ and $\Gamma(\mathfrak{J}_1)$ is the restriction of $\Gamma(\mathfrak{J})$ on $(\mathfrak{J}_1)_{ss}.$
\end{Lem}

Let $\mathfrak{J}$ be a non-nilpotent finite-dimensional Jordan algebra. Then $\mathfrak{J}$ has a nonzero idempotent $e$, and $\mathfrak{J}$ admits the following decomposition:
$$\mathfrak{J} = \mathfrak{J}_0 \oplus \mathfrak{J}_1 \oplus \mathfrak{J}_2,$$
where $\mathfrak{J}_i = \{x \in \mathfrak{J}: ex = \frac{i}{2}x\},$ called the \textit{Peirce decomposition of $\mathfrak{J}$ with respect to $e$.} 
The spaces $\mathfrak{J}_0, \mathfrak{J}_1$ and $\mathfrak{J}_2$ satisfy the following multiplicative relations:
\begin{equation}
\label{pd}
\mathfrak{J}_0\mathfrak{J}_0 \subseteq \mathfrak{J}_0, \mathfrak{J}_2\mathfrak{J}_2 \subseteq \mathfrak{J}_2, \mathfrak{J}_0\mathfrak{J}_2 = 0, \mathfrak{J}_0\mathfrak{J}_1 \subseteq \mathfrak{J}_1, \mathfrak{J}_2\mathfrak{J}_1 \subseteq \mathfrak{J}_1, \mathfrak{J}_1\mathfrak{J}_1 \subseteq \mathfrak{J}_0 + \mathfrak{J}_2.
\end{equation}
Analogously, if $e_1 + \ldots + e_n =  1\in \mathfrak{J}$ is the sum of $n$ orthogonal idempotents, then $\mathfrak{J}$ has the following Peirce decomposition with respect to $e_1, \ldots, e_n:$
$$\mathfrak{J} = \bigoplus_{i,j = 1}^n \mathfrak{J}_{ij},$$
where
$$\mathfrak{J}_{ii} = \{x \in \mathfrak{J}: e_ix = 0, e_jx = 0, j \neq i\},$$
$$\mathfrak{J}_{ij} = \{x \in \mathfrak{J}: e_ix = e_jx = \frac{1}{2}x, e_kx = 0, k \neq i, j\} = \mathfrak{J}_{ji}, i \neq j.$$

As a corollary of the previous lemma, we get 
\begin{Lem}[\cite{K1}]
\label{pd_preserve}
If $\mathfrak{J} \to \mathfrak{J}_1$ for $\mathfrak{J}_1$ non-nilpotent, then $(\mathfrak{J}_1)_{ss}$ is a subalgebra of $\mathfrak{J}_{ss}$ and the degeneration preserves the corresponding Peirce decomposition.
\end{Lem}




\section{Degenerations of low dimensional Jordan algebras}\label{deg}
In this section we apply the methods described in the previous sections to completely describe the degenerations of low-dimensional Jordan algebras.

\subsection{$\mathfrak{Jord}_2$}
The description of degenerations in $\mathfrak{Jord}_2$ follows from \cite{kv16}.
For a more convenient notation, we will list here all $2$-dimensional Jordan algebras:

\begin{center}
Table 1. {\it Jordan algebras of dimension $2$.}


\footnotesize

$$\begin{array}{|l|l|c|} 
\hline
\mathbb{A}  & \mbox{multiplication tables} & \mathfrak{Der} \\ 

\hline \hline 

\mathfrak{B}_1^{A}& \begin{array}{ll}e_1^2 = e_1, & e_1n_1 =n_1 \end{array} & 1\\ 
\hline

\mathfrak{B}_2 & \begin{array}{ll}e_1^2 = e_1, & e_1n_1 =\frac{1}{2}n_1\end{array} & 2\\ 

\hline
\mathfrak{B}_3^{AN}& \begin{array}{ll}n_1^2 = n_2\end{array} & 2\\ 

\hline
\mathfrak{B}_4^{A}   & \begin{array}{ll}e_1^2 = e_1, & e_2^2 = e_2 \end{array} & 0\\ 

\hline
\mathfrak{B}_5^{A}  & \begin{array}{ll}e_1^2 = e_1 \end{array} & 1\\ 

\hline

\end{array}$$
\end{center}

\normalsize

\begin{Th}
\label{2d}
The graph of primary degenerations for $\mathfrak{Jord}_2$  has the following form: 
\end{Th}

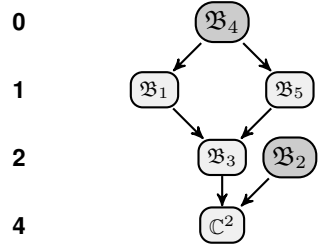
\footnotesize

$${\begin{tikzpicture}[->,>=stealth',shorten >=0.0cm,auto,node distance=0.9cm,
                    thick,main node/.style={rectangle,draw,fill=gray!12,rounded corners=1ex,font=\sffamily \tiny 
                    \bfseries },rigid node/.style={rectangle,draw,fill=black!20,rounded corners=1.5ex,font=\sffamily \bf \bfseries },style={draw,font=\sffamily \scriptsize \bfseries }]
\node (0)   {0};

\node (00a1) [right of=0] {};
\node (00a2) [right of=00a1] {};
\node (00a3) [right of=00a2] {};
\node (00a4) [right of=00a3] {};

\node (01) [below of=0] {1};

\node (01a1) [right of=01] {};
\node (01a2) [right of=01a1] {};
\node (01a3) [right of=01a2] {};
\node (01a4) [right of=01a3] {};

\node (02) [below of=01] {2};

\node (02a1) [right of=02] {};
\node (02a2) [right of=02a1] {};
\node (02a3) [right of=02a2] {};
\node (02a4) [right of=02a3] {};

\node (04) [below of=02] {4};

\node (04a1) [right of=04] {};
\node (04a2) [right of=04a1] {};
\node (04a3) [right of=04a2] {};
\node (04a4) [right of=04a3] {};

\node  [rigid node] (b4) [right of=00a2] {$\mathfrak{B}_{4}$};
\node  [main node] (b1) [right of=01a1] {$\mathfrak{B}_{1}$};	
\node  [rigid node] (b2) [right of=02a3] {$\mathfrak{B}_{2}$};	
\node  [main node] (b5) [right of=01a3] {$\mathfrak{B}_{5}$};
\node  [main node] (b3) [right of=02a2] {$\mathfrak{B}_{3}$};
\node  [main node] (c2) [right of=04a2] {$\mathbb{C}^{2}$};

\path[every node/.style={font=\sffamily\small}]

(b4) edge   node[above] {} (b1)
(b4) edge   node[above] {} (b5)
(b5) edge   node[above] {} (b3)
(b1) edge   node[above] {} (b3)
(b2) edge   node[above] {} (c2)
(b3) edge   node[above] {} (c2);

\end{tikzpicture}}$$

\normalsize

\begin{corollary}\label{ir_j2} The irreducible components  of
$\mathfrak{Jord}_2$ are
$$
\begin{aligned}
\mathcal{C}_1   &=\overline{ O(\mathfrak{B}_{2})  }=  \{ \mathfrak{B}_{2}, \mathbb{C}^2 \};\\
\mathcal{C}_2   &=\overline{ O(\mathfrak{B}_{4})  }=  \{ \mathfrak{B}_{1}, \mathfrak{B}_{3}, \mathfrak{B}_{4}, \mathfrak{B}_{5}, \mathbb{C}^2\}.\\
\end{aligned}
$$ 
In particular, 
$Rig(\mathfrak{Jord}_2)= 
\{ \mathfrak{B}_{2}, \mathfrak{B}_{4}  \}.$
\end{corollary}

\subsection{$\mathfrak{Jord}_3$}

In the Table 2 we list all  algebras from $\mathfrak{Jord}_3.$

\footnotesize

\begin{center}
Table 2. {\it Jordan algebras of dimension $3$.}
\begin{equation*}
\begin{array}{|l|l|l|l|c|c|} 
\hline
\mathbb{A}  &  \mbox{\cite{KS07}}  &  \mbox{ Decomposition } & \mbox{ multiplication tables } & \mathfrak{Der} & \mathfrak{Rad} 
\\ 

\hline \hline 
\mathbb{T}_{01}^{AUS}    &A_{11}& \mathbb{C}e_1 \oplus \mathbb{C}e_2 \oplus \mathbb{C}e_3  & \begin{array}{c} e_1^2 = e_1, e_2^2 = e_2, e_3^2 = e_3  \end{array}
& 0 & 0 \\

\hline
\mathbb{T}_{02}^{US}& J_1&-& \begin{array}{c} e_1^2 = e_1, e_2^2 =e_2, e_3^2 = e_1 + e_2 , e_1e_3 = \frac{1}{2}e_3, e_2e_3 = \frac{1}{2}e_3\end{array}& 1 & 0\\ 

\hline 
\hline
\mathbb{T}_{03}^{AU} & A_{12}& \mathfrak{B}_1 \oplus \mathbb{C}e_2 & \begin{array}{c} e_1^2 = e_1,   e_2^2 = e_2, e_1n_1 = n_1 \end{array} & 1 & 1\\

\hline
\mathbb{T}_{04}^{U} &J_2&- & \begin{array}{c}   e_1^2= e_1, e_2^2=e_2, e_1n_1 = \frac{1}{2}n_1, e_2n_1 = \frac{1}{2}n_1  \end{array}
& 2 & 1 \\  

\hline

\mathbb{T}_{05} & J_3&  \mathfrak{B}_2 \oplus \mathbb{C}e_1 & \begin{array}{c}e_1^2 = e_1, e_2^2 = e_2,  e_1n_1 = \frac{1}{2}n_1  \end{array} & 2 & 1\\

\hline
\mathbb{T}_{06}^{A}
&A_1& \mathbb{C}e_1 \oplus \mathbb{C}e_2\oplus \mathbb{C}n_1 & \begin{array}{c} e_1^2 = e_1, e_2^2 = e_2 \end{array} & 1
& 1 \\ 
\hline
\hline

\mathbb{T}_{07}^{AU}& A_{13} & - &\begin{array}{c}e_1^2 = e_1, e_1n_1 =n_1, e_1n_2=n_2, n_1^2 = n_2 \end{array} & 2 & 2\\ 
\hline
\mathbb{T}_{08}^{AU} & A_{14}& - & \begin{array}{c}e_1^2 = e_1, e_1n_1 =n_1, e_1n_2=n_2\end{array}& 4  & 2\\ 
\hline

\mathbb{T}_{09}^{A} & A_2& \mathfrak{B}_1 \oplus \mathbb{C}n_2 & \begin{array}{c} e_1^2 = e_1, e_1n_1 = n_1\end{array} & 2 & 2\\ 
\hline

\mathbb{T}_{10} & J_7&- & \begin{array}{c} e_1^2 = e_1, e_1n_1 = \frac{1}{2}n_1, e_1n_2 = n_2, n_1^2 = n_2   \end{array}
& 2 & 2 \\ 
\hline
\mathbb{T}_{11}& J_4& -& \begin{array}{c}e_1^2 = e_1,e_1n_1 =\frac{1}{2}n_1,e_1n_2 = n_2\end{array}& 3& 2 \\ 
\hline
\mathbb{T}_{12}& J_5&- & \begin{array}{c}e_1^2 = e_1, e_1n_1 = \frac{1}{2}n_1, e_1n_2 = \frac{1}{2}n_2 \end{array} 
& 6 & 2\\

\hline
\mathbb{T}_{13} & J_6& -& \begin{array}{c}e_1^2 = e_1, e_1n_1 = \frac{1}{2}n_1, n_1^2 = n_2  \end{array}& 2 & 2 \\ 
\hline

\mathbb{T}_{14} & J_8& \mathfrak{B}_2 \oplus \mathbb{C}n_2 &  \begin{array}{c} e_1^2 = e_1, e_1n_1 = \frac{1}{2}n_1 \end{array}& 
3 & 2\\

\hline
\mathbb{T}_{15}^{A} & A_3& \mathfrak{B}_3 \oplus \mathbb{C}e_1  & \begin{array}{c} e_1^2 = e_1, n_1^2 = n_2\end{array}& 2 & 2\\
\hline

\mathbb{T}_{16}^{A} & A_5& \mathbb{C}e_1 \oplus \mathbb{C}n_1 \oplus \mathbb{C}n_2   & \begin{array}{c}e_1^2 = e_1  \end{array} & 4 & 2\\ 
\hline
\hline 
\mathbb{T}_{17}^{AN} & A_4& -& \begin{array}{c}n_1^2 = n_2, n_1n_2 = n_3\end{array}& 3 & 3\\ 

\hline

\mathbb{T}_{18}^{AN}& A_6&-& \begin{array}{c}
n_1n_2 = n_3\end{array} & 4 & 3\\ 

\hline
\mathbb{T}_{19}^{AN} & A_7& \mathfrak{B}_3 \oplus \mathbb{C}n_1 &\begin{array}{c} n_1^2 = n_2 \end{array} & 5  & 3 \\ \hline

\end{array}
\end{equation*}
\end{center}

\normalsize 

\begin{Th}
\label{3d}
The graph of primary degenerations for $\mathfrak{Jord}_3$  has the following form: 
\end{Th}
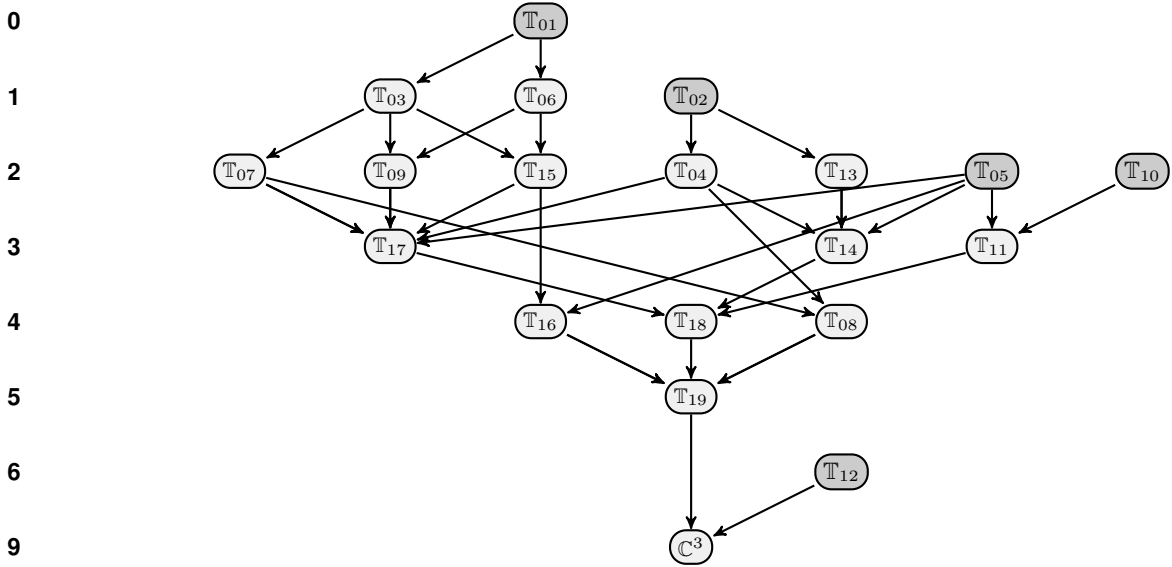
\scriptsize 

\begin{center}
\begin{tikzpicture}[->,>=stealth',shorten >=0.0cm,auto,node distance=1cm,thick,
                    main node/.style={rectangle,draw,fill=gray!12,rounded corners=1.5ex,font=\sffamily \tiny \bfseries },
                    rigid node/.style={rectangle,draw,fill=black!20,rounded corners=1.5ex,font=\sffamily \bf \bfseries },
                    style={draw,font=\sffamily \scriptsize \bfseries }]
\node (0)   {0};

\node (00a1) [right of=0] {};
\node (00a2) [right of=00a1] {};
\node (00a3) [right of=00a2] {};
\node (00a4) [right of=00a3] {};
\node (00a5) [right of=00a4] {};
\node (00a6) [right of=00a5] {};
\node (00a7) [right of=00a6] {};
\node (00a8) [right of=00a7] {};
\node (00a9) [right of=00a8] {};
\node (00a10) [right of=00a9] {};
\node (00a11) [right of=00a10] {};
\node (00a12) [right of=00a11] {};
\node (00a13) [right of=00a12] {};
\node (00a14) [right of=00a13] {};
\node (00a15) [right of=00a14] {};
\node (00a16) [right of=00a15] {};
\node (00a17) [right of=00a16] {};
\node (00a18) [right of=00a17] {};
\node (00a19) [right of=00a18] {};
\node (00a20) [right of=00a19] {};
\node (00a21) [right of=00a20] {};
\node (00a22) [right of=00a21] {};
\node (00a23) [right of=00a22] {};
\node (00a24) [right of=00a23] {};
\node (00a25) [right of=00a24] {};
\node (00a26) [right of=00a25] {};
\node (00a27) [right of=00a26] {};

\node (01) [below of=0] {1};

\node (01a1) [right of=01] {};
\node (01a2) [right of=01a1] {};
\node (01a3) [right of=01a2] {};
\node (01a4) [right of=01a3] {};
\node (01a5) [right of=01a4] {};
\node (01a6) [right of=01a5] {};
\node (01a7) [right of=01a6] {};
\node (01a8) [right of=01a7] {};
\node (01a9) [right of=01a8] {};
\node (01a10) [right of=01a9] {};
\node (01a11) [right of=01a10] {};
\node (01a12) [right of=01a11] {};
\node (01a13) [right of=01a12] {};
\node (01a14) [right of=01a13] {};
\node (01a15) [right of=01a14] {};

\node (02) [below of=01] {2};

\node (02a1) [right of=02] {};
\node (02a2) [right of=02a1] {};
\node (02a3) [right of=02a2] {};
\node (02a4) [right of=02a3] {};
\node (02a5) [right of=02a4] {};
\node (02a6) [right of=02a5] {};
\node (02a7) [right of=02a6] {};
\node (02a8) [right of=02a7] {};
\node (02a9) [right of=02a8] {};
\node (02a10) [right of=02a9] {};
\node (02a11) [right of=02a10] {};
\node (02a12) [right of=02a11] {};
\node (02a13) [right of=02a12] {};
\node (02a14) [right of=02a13] {};
\node (02a15) [right of=02a14] {};
\node (02a16) [right of=02a15] {};
\node (02a17) [right of=02a16] {};
\node (02a18) [right of=02a17] {};
\node (02a19) [right of=02a18] {};
\node (02a20) [right of=02a19] {};
\node (02a21) [right of=02a20] {};
\node (02a22) [right of=02a21] {};
\node (02a23) [right of=02a22] {};
\node (02a24) [right of=02a23] {};
\node (02a25) [right of=02a24] {};
\node (02a26) [right of=02a25] {};
\node (02a27) [right of=02a26] {};

\node (03)[below of=02]{3};

\node (03a1) [right of=03] {};
\node (03a2) [right of=03a1] {};
\node (03a3) [right of=03a2] {};
\node (03a4) [right of=03a3] {};
\node (03a5) [right of=03a4] {};
\node (03a6) [right of=03a5] {};
\node (03a7) [right of=03a6] {};
\node (03a8) [right of=03a7] {};
\node (03a9) [right of=03a8] {};
\node (03a10) [right of=03a9] {};
\node (03a11) [right of=03a10] {};
\node (03a12) [right of=03a11] {};
\node (03a13) [right of=03a12] {};
\node (03a14) [right of=03a13] {};
\node (03a15) [right of=03a14] {};
\node (03a16) [right of=03a15] {};
\node (03a17) [right of=03a16] {};
\node (03a18) [right of=03a17] {};
\node (03a19) [right of=03a18] {};
\node (03a20) [right of=03a19] {};
\node (03a21) [right of=03a20] {};
\node (03a22) [right of=03a21] {};
\node (03a23) [right of=03a22] {};
\node (03a24) [right of=03a23] {};
\node (03a25) [right of=03a24] {};
\node (03a26) [right of=03a25] {};
\node (03a27) [right of=03a26] {};
*\node (03a28) [right of=03a27] {};

\node (04) [below of=03] {4};

\node (04a1) [right of=04] {};
\node (04a2) [right of=04a1] {};
\node (04a3) [right of=04a2] {};
\node (04a4) [right of=04a3] {};
\node (04a5) [right of=04a4] {};
\node (04a6) [right of=04a5] {};+
\node (04a7) [right of=04a6] {};
\node (04a8) [right of=04a7] {};
\node (04a9) [right of=04a8] {};
\node (04a10) [right of=04a9] {};
\node (04a11) [right of=04a10] {};
\node (04a12) [right of=04a11] {};
\node (04a13) [right of=04a12] {};
\node (04a14) [right of=04a13] {};
\node (04a15) [right of=04a14] {};
\node (04a16) [right of=04a15] {};
\node (04a17) [right of=04a16] {};
\node (04a18) [right of=04a17] {};
\node (04a19) [right of=04a18] {};
\node (04a20) [right of=04a19] {};
\node (04a21) [right of=04a20] {};
\node (04a22) [right of=04a21] {};
\node (04a23) [right of=04a22] {};
\node (04a24) [right of=04a23] {};
\node (04a25) [right of=04a24] {};
\node (04a26) [right of=04a25] {};
\node (04a27) [right of=04a26] {};
\node (04a28) [right of=04a27] {};

\node (05) [below of=04] {5};

\node (05a1) [right of=05] {};
\node (05a2) [right of=05a1] {};
\node (05a3) [right of=05a2] {};
\node (05a4) [right of=05a3] {};
\node (05a5) [right of=05a4] {};
\node (05a6) [right of=05a5] {};
\node (05a7) [right of=05a6] {};
\node (05a8) [right of=05a7] {};
\node (05a9) [right of=05a8] {};
\node (05a10) [right of=05a9] {};
\node (05a11) [right of=05a10] {};
\node (05a12) [right of=05a11] {};
\node (05a13) [right of=05a12] {};
\node (05a14) [right of=05a13] {};
\node (05a15) [right of=05a14] {};
\node (05a16) [right of=05a15] {};
\node (05a17) [right of=05a16] {};
\node (05a18) [right of=05a17] {};
\node (05a19) [right of=05a18] {};
\node (05a20) [right of=05a19] {};
\node (05a21) [right of=05a20] {};
\node (05a22) [right of=05a21] {};
\node (05a23) [right of=05a22] {};
\node (05a24) [right of=05a23] {};
\node (05a25) [right of=05a24] {};
\node (05a26) [right of=05a25] {};
\node (05a27) [right of=05a26] {};
\node (05a28) [right of=05a27] {};

\node (06)[below of=05] {6};

\node (09) [below of=06] {9};

\node (06a1) [right of=06] {};
\node (06a2) [right of=06a1] {};
\node (06a3) [right of=06a2] {};
\node (06a4) [right of=06a3] {};
\node (06a5) [right of=06a4] {};
\node (06a6) [right of=06a5] {};
\node (06a7) [right of=06a6] {};
\node (06a8) [right of=06a7] {};
\node (06a9) [right of=06a8] {};
\node (06a10) [right of=06a9] {};
\node (06a11) [right of=06a10] {};
\node (06a12) [right of=06a11] {};
\node (06a13) [right of=06a12] {};
\node (06a14) [right of=06a13] {};
\node (06a15) [right of=06a14] {};
\node (06a16) [right of=06a15] {};
\node (06a17) [right of=06a16] {};
\node (06a18) [right of=06a17] {};
\node (06a19) [right of=06a18] {};
\node (06a20) [right of=06a19] {};
\node (06a21) [right of=06a20] {};
\node (06a22) [right of=06a21] {};
\node (06a23) [right of=06a22] {};
\node (06a24) [right of=06a23] {};
\node (06a25) [right of=06a24] {};
\node (06a26) [right of=06a25] {};
\node (06a27) [right of=06a26] {};
\node (06a28) [right of=06a27] {};

\node  [main node] (t1) [right of = 02a2] {$\mathbb{T}_{07}$};

\node  [main node] (t12) [right of=01a4] {$\mathbb{T}_{03}$};	
\node  [main node] (t15) [right of=02a4] {$\mathbb{T}_{09}$};	
\node  [main node] (t3)  [right of=03a4] {$\mathbb{T}_{17}$};

\node  [rigid node] (t11) [right of=00a6] {$\mathbb{T}_{01}$};	
\node  [main node] (t13) [right of=01a6] {$\mathbb{T}_{06}$};	
\node  [main node] (t16) [right of=02a6] {$\mathbb{T}_{15}$};
\node  [main node] (t17) [right of=04a6] {$\mathbb{T}_{16}$};

\node  [rigid node] (t5)  [right of=01a8] {$\mathbb{T}_{02}$};		
\node  [main node] (t10) [right of=02a8] {$\mathbb{T}_{04}$};
\node  [main node] (t4)  [right of=04a8] {$\mathbb{T}_{18}$};
\node  [main node] (t19) [right of=05a8] {$\mathbb{T}_{19}$};

\node  [main node] (t8)  [right of=02a10] {$\mathbb{T}_{13}$};
\node  [main node] (t18) [right of=03a10] {$\mathbb{T}_{14}$};
\node  [main node] (t2)  [right of=04a10] {$\mathbb{T}_{08}$};
\node  [rigid node] (t7) [right of=06a10] {$\mathbb{T}_{12}$};

\node  [rigid node] (t14) [right of=02a12] {$\mathbb{T}_{05}$};
\node  [main node] (t6)  [right of=03a12] {$\mathbb{T}_{11}$};

\node  [rigid node] (t9) [right of=02a14] {$\mathbb{T}_{10}$};

\node  [main node] (t20) [below of=06a9] {$\mathbb{C}^3$};

\path[every node/.style={font=\sffamily\small}]

(t1) edge   node[above] {} (t2)
(t1) edge   node[above] {} (t3)
(t1) edge   node[above] {} (t3)

(t2) edge   node[above] {} (t19)
(t2) edge   node[above] {} (t19)

(t3) edge   node[above] {} (t4)

(t4) edge   node[above] {} (t19)

(t5) edge   node[above] {} (t8)
(t5) edge   node[above] {} (t10)

(t6) edge   node[above] {} (t4)

(t7) edge   node[above] {} (t20)

(t8) edge   node[above] {} (t18)
(t8) edge   node[above] {} (t18)

(t9) edge   node[above] {} (t6)

(t10) edge   node[above] {} (t18)
(t10) edge   node[above] {} (t2)

(t11) edge   node[above] {} (t12)
(t11) edge   node[above] {} (t13)

(t12) edge   node[above] {} (t1)
(t12) edge   node[above] {} (t15)
(t12) edge   node[above] {} (t16)

(t13) edge   node[above] {} (t15)
(t13) edge   node[above] {} (t16)

(t14) edge   node[above] {} (t6)
(t14) edge   node[above] {} (t17)
(t14) edge   node[above] {} (t18)

(t15) edge   node[above] {} (t3)
(t15) edge   node[above] {} (t3)

(t16) edge   node[above] {} (t3)
(t16) edge   node[above] {} (t17)

(t17) edge   node[above] {} (t19)
(t17) edge   node[above] {} (t19)

(t18) edge   node[above] {} (t4)

(t19) edge   node[above] {} (t20);

\end{tikzpicture}

\end{center}
\normalsize

 \begin{Proof} 
First we note that the set of rigid 3-dimensional Jordan algebras is $\{\mathbb{T}_{01}, \mathbb{T}_{02},\mathbb{T}_{05} ,\mathbb{T}_{10}, \mathbb{T}_{12}\}$ (follows from \cite{KS07}).

Then we list all primary degenerations of $\mathfrak{Jord}_3:$
\begin{enumerate}
\item Primary degenerations from  \cite{KS07}:
\begin{eqnarray*}
\mathbb{T}_{01} \to \mathbb{T}_{03}, \mathbb{T}_{01} \to \mathbb{T}_{06}, \mathbb{T}_{02} \to \mathbb{T}_{04}, \mathbb{T}_{02} \to \mathbb{T}_{13},  \mathbb{T}_{03} \to \mathbb{T}_{07},  \mathbb{T}_{04} \to \mathbb{T}_{14}, \mathbb{T}_{06} \to \mathbb{T}_{09},\\ \mathbb{T}_{06} \to \mathbb{T}_{15}, \mathbb{T}_{07} \to \mathbb{T}_{08},  \mathbb{T}_{10} \to \mathbb{T}_{11},  \mathbb{T}_{15} \to \mathbb{T}_{16},  \mathbb{T}_{15} \to \mathbb{T}_{17}, \mathbb{T}_{17} \to \mathbb{T}_{18}, \mathbb{T}_{18} \to \mathbb{T}_{19}. 
\end{eqnarray*}

\item The list of all possible degenerations is completed in the table below.

\end{enumerate}

\begin{center} \footnotesize Table. {\it Primary degenerations of 3-dimensional Jordan algebras.}
$$\begin{array}{|l|lll|}
\hline
\mbox{Degenerations}  &  \multicolumn{3}{c|}{\mbox{Parametrized bases}}  \\
\hline
\hline

\mathbb{T}_{03} \to \mathbb{T}_{09} & E_1^t = e_1 &E_2^t = n_1 &E_3^t = te_2 \\
\hline

\mathbb{T}_{03} \to \mathbb{T}_{15} & E_1^t = e_2 & E_2^t = te_1 + n_1 &E_3^t = -t^2n_1\\
\hline

\mathbb{T}_{04} \to \mathbb{T}_{08} & E_1^t = e_1+e_2 & E_2^t = t(e_1-e_2) &E_3^t = tn_1 \\
\hline


\mathbb{T}_{05} \to \mathbb{T}_{11} & E_1^t = e_1 + e_2 & E_2^t = n_1 & E_3^t = te_1 \\
\hline

\mathbb{T}_{05} \to \mathbb{T}_{14} & E_1^t = e_1 & E_2^t = n_1 & E_3^t = te_2 \\
\hline

\mathbb{T}_{05} \to \mathbb{T}_{16} & E_1^t = e_2 & E_2^t = te_1 &E_3^t =  n_1\\
\hline


\mathbb{T}_{07} \to \mathbb{T}_{17} & E_1^t = te_1 + n_1  & E_2^t = tn_1 + n_2  & E_3^t = tn_2 \\
\hline

\mathbb{T}_{08} \to \mathbb{T}_{19} &  E_1^t = te_1 + n_1 & E_2^t = tn_1, & E_3^t = n_2  \\
\hline

\mathbb{T}_{09} \to \mathbb{T}_{17} & E_1^t = te_1 + n_1 + n_2 & E_2^t = tn_1 - tn_2 & E_3^t = t^2n_1 \\
\hline

\mathbb{T}_{11} \to \mathbb{T}_{18} & E_1^t = te_1 & E_2^t = n_1 + 2n_2 & E_3^t = tn_2\\
\hline

\mathbb{T}_{13} \to \mathbb{T}_{14} & E_1^t = e_1 & E_2^t = tn_1 & E_3^t = n_2\\
\hline

\mathbb{T}_{14} \to \mathbb{T}_{18} & E_1^t = te_1  & E_2^t = n_1 - 2n_2 & E_3^t = tn_2\\
\hline

\mathbb{T}_{16} \to \mathbb{T}_{19} & E_1^t = te_1+n_1  & E_2^t = tn_1 & E_3^t = n_2 \\
\hline

\hline
\end{array}$$
\end{center}

Finally, we give a descriprion of all primary non-degenerations.
This description  has four  parts:
\begin{enumerate}
    \item Non-degenerations which follows from lemma \ref{pd_preserve}: 
\begin{eqnarray*}
    \mathbb{T}_{02} \not\to \mathbb{T}_{09}, \mathbb{T}_{11}, \mathbb{T}_{16}; \mathbb{T}_{05} \not\to \mathbb{T}_{08}; \mathbb{T}_{06} \not\to \mathbb{T}_{07}, \mathbb{T}_{08}; \mathbb{T}_{07} \not\to \mathbb{T}_{16}; \\ 
    \mathbb{T}_{09} \not\to \mathbb{T}_{16}; \mathbb{T}_{10} \not\to \mathbb{T}_{08}, \mathbb{T}_{14}, \mathbb{T}_{16}; \mathbb{T}_{13} \not\to \mathbb{T}_{08}. 
\end{eqnarray*}

\item To prove the assertions $\mathbb{T}_{10} \not\to \mathbb{T}_{17}$ and $ \mathbb{T}_{13} \not\to \mathbb{T}_{17} $ let us consider the set
$$
\mathcal{R}=\left\{\mu\in\mathfrak{Jord}_3\left|\begin{array}{c} c_{11}^2=0, S_1S_3+ S_2^2 \subseteq S_3, S_2S_3=0 \end{array}\right.\right\}.
$$
It is not difficult to show that $\mathcal{R}$ is a closed subset of $\mathfrak{Jord}_3$ stable with respect to the action of the subgroup of upper triangular matrices and that $\mathcal{R}$ contains the structures $\mathbb{T}_{10}$ and $\mathbb{T}_{13}$ (to see this, it is enough to consider the elements of $\mathfrak{Jord}_3$ corresponding to these algebras in the bases $e_1$, $n_1$, $n_2$).
It is also not difficult to show that $\mathcal{R}\cap O(\mathbb{T}_{17})= \emptyset.$ 

\item 

To prove the assertion $\mathbb{T}_{02} \not\to \mathbb{T}_{07}$ let us consider the set
$$
\mathcal{R}=\left\{\mu\in\mathfrak{Jord}_3\left|\begin{array}{c} c_{11}^2=c_{23}^3=c_{13}^3=c_{13}^2=0,  S_2^2 \subseteq S_2, S_3^2 \subseteq S_3 \end{array}\right.\right\}.
$$
It is not difficult to show that $\mathcal{R}$ is a closed subset of $\mathfrak{Jord}_3$  stable with respect to the action of the subgroup of upper triangular matrices and that $\mathcal{R}$ contains the structure $\mathbb{T}_{02}.$  
To see this it is enough to consider the element corresponding to this algebra in the basis $e_3$, $e_2$, $e_1+e_2.$
It is also not difficult to show that $\mathcal{R}\cap O(\mathbb{T}_{07})= \emptyset.$

\item 

To prove the assertion $\mathbb{T}_{04}, \mathbb{T}_{05} \not\to \mathbb{T}_{17}$ let us consider the set
$$
\mathcal{R}=\left\{\mu\in\mathfrak{Jord}_3\left|\begin{array}{c} c_{11}^2c_{12}^3=(\frac{c_{11}^2}{2}-c_{12}^2) c_{11}^3,  
S_1S_2 \subseteq S_2, S_3^2=0, S_1S_3 \subseteq S_3 \end{array}\right.\right\}.
$$
It is not difficult to show that $\mathcal{R}$ is a closed subset of $\mathfrak{Jord}_3$  stable with respect to the action of the subgroup of upper triangular matrices and that $\mathcal{R}$ contains the structure $\mathbb{T}_{04}, \mathbb{T}_{05}.$  
To see this it is enough to consider the element corresponding to this algebra in the basis $e_1$, $e_2$, $n_1.$
It is also not difficult to show that $\mathcal{R}\cap O(\mathbb{T}_{17})= \emptyset.$

\end{enumerate}

\end{Proof}

\begin{corollary}\label{ir_j3} The irreducible components  of
$\mathfrak{Jord}_3$ are
$$
\begin{aligned}
\mathcal{C}_1   &=\overline{ O(\mathbb{T}_{01})  }=  \{ \mathbb{T}_{01}, \mathbb{T}_{03}, \mathbb{T}_{06}, \mathbb{T}_{07}, \mathbb{T}_{08}, \mathbb{T}_{15}, \mathbb{T}_{16}, \mathbb{T}_{17}, \mathbb{T}_{18},  \mathbb{T}_{09},  \mathbb{T}_{19}, \mathbb{C}^3 \};\\
\mathcal{C}_2   &=\overline{ O(\mathbb{T}_{02})  }=  \{ \mathbb{T}_{02}, \mathbb{T}_{04}, \mathbb{T}_{08}, \mathbb{T}_{13}, \mathbb{T}_{14},   \mathbb{T}_{18}, \mathbb{T}_{19}, \mathbb{C}^3 \};\\
\mathcal{C}_3   &=\overline{ O(\mathbb{T}_{05})  }=  \{ \mathbb{T}_{05}, \mathbb{T}_{11}, \mathbb{T}_{14}, \mathbb{T}_{16},  \mathbb{T}_{18}, \mathbb{T}_{19}, \mathbb{C}^3\};\\
\mathcal{C}_4   &=\overline{ O(\mathbb{T}_{10})  }=  \{ \mathbb{T}_{10}, \mathbb{T}_{11}, \mathbb{T}_{18}, \mathbb{T}_{19}, \mathbb{C}^3\};\\
\mathcal{C}_5   &=\overline{ O(\mathbb{T}_{12})  }=  \{ \mathbb{T}_{12}, \mathbb{C}^3\}.\\
\end{aligned}
$$ 
In particular, 
$Rig(\mathfrak{Jord}_3)= 
\{\mathbb{T}_{01}, \mathbb{T}_{02},\mathbb{T}_{05} ,\mathbb{T}_{10}, \mathbb{T}_{12}\}.$
\end{corollary}

\section{''Marginal'' algebras}\label{marginal}

\subsection{Jordan algebras}\label{mar} In this section we consider a special series of Jordan algebras $\mathfrak{J}_n$ of dimension $n, n \geq 2.$ The algebra $\mathfrak{J}_n$ has a basis $e, n_1, \ldots n_{n-1},$ and its multiplication table is given by:
$$e^2 = e, en_i = \frac{1}{2}n_i, n_in_j = 0, \text{ where } i, j = 1,\ldots,n-1.$$
In the $2$-dimensional case we denoted this algebra as $\mathfrak{B}_2$, and in the $3$-dimensional as $\mathbb{T}_{12}$. These algebras have a certain ''marginal'' behavior in the varieties $\mathfrak{Jord}_n, n = 2, 3$: the algebras $\mathfrak{B}_2$ and $\mathbb{T}_{12}$ are rigid, but only degenerate to the algebra with zero multiplication. The goal of this section is to prove that the algebras $\mathfrak{J}_n$ exhibit this behavior for all $n \geq 2$.

The Peirce decomposition of this algebra 
with respect to the idempotent $e$ is the following:
$$(\mathfrak{J}_n)_0 = 0, (\mathfrak{J}_n)_1 = \langle n_1, \ldots, n_{n-1}\rangle, (\mathfrak{J}_n)_2 = \langle e \rangle.$$
Firstly, we show that the algebra $\mathfrak{J}_n$ is of level one, that is, if $\mathfrak{J}_n \to \mathfrak{J}'$, where $\mathfrak{J}' \in \mathfrak{Jord}_n,$ then $\mathfrak{J}'$ is an algebra with zero multiplication. Note that this was already proved in the paper \cite{khud13}, but here we present a simpler proof of this fact.

First of all, we will need the following technical statement:
\begin{Lem}
\label{marg_der}
The dimension of the derivation algebra of $\mathfrak{J}_n$ is equal to $n^2 - n.$
\end{Lem}
\begin{Proof} Let $d$ be a derivation of $\mathfrak{J}_n.$ Then
$$d(e) = d(e^2) = 2ed(e),$$
therefore $d(e) \in (\mathfrak{J}_n)_1 = \langle n_1, \ldots n_{n-1} \rangle.$ Now let $a \in (\mathfrak{J}_n)_1.$ Then
$$\frac{1}{2}d(a) = d(ea) = d(e)a + ed(a) = ed(a),$$
since $d(e) \in (\mathfrak{J}_n)_1$ and $(\mathfrak{J}_n)_1^2 = 0$. The defining relation for $d$ to be a derivation is therefore easily seen to hold for $a, b \in (\mathfrak{J}_n)_1.$ Thus a linear map $d$ on $\mathfrak{J}_n$ is a derivation if and only if its image is contained in $(\mathfrak{J}_n)_1.$ Hence $\dim \mathfrak{Der}(\mathfrak{J}_n) = n(n-1).$ 
\end{Proof}

Alternatively, we could use the proof of this lemma in the work \cite{Kap} ($\mathfrak{J}_n$ is a ''$T$-algebra'' in the terminology of this paper).

Now we can show the following statement:
\begin{Lem}
\label{marg_deg}
Let $\mathfrak{J}_n \to \mathfrak{J}'$ be a proper degeneration, where $\mathfrak{J}'$ is an algebra of dimension $n.$ Then $\mathfrak{J}'$ has zero multiplication. In other words, $\mathfrak{J}_n$ is of level one.
\end{Lem}
\begin{Proof} From the results of the paper \cite{Kap}, an algebra of dimension $n$ has a derivation algebra of dimension less or equal than $ n^2 - n.$ Therefore the previous lemma and the fact that the dimension of derivation algebra grows strictly under proper degeneration imply the statement of the lemma.
\end{Proof}

\medskip
Another peculiar property of the algebra $\mathfrak{J}_n$ is that it is rigid in the variety $\mathfrak{J}_n$. We will prove this fact by showing that the second cohomology 
group of $\mathfrak{J}_n$ with coefficients in itself is zero. Let us recall the definition of the second cohomology group of a Jordan algebra: Let $\mathfrak{J}$ be a Jordan algebra. Let $\mathbb{V}$ be the space of the bilinear maps $h: \mathfrak{J}\times\mathfrak{J} \to \mathfrak{J}$ such that (see \cite{Jac})
\begin{equation}
\label{cohom_1}
h(a,b) = h(b,a),
\end{equation}
\begin{equation}
\label{cohom_2}
(h(a,a)b) a + h(a^2,b)a + h(a^2b,a) + a^2h(b,a) + h(a,a)(ba) + h(a^2,ba) = 0.
\end{equation}
For example, if $\mu: \mathfrak{J} \to \mathfrak{J}$ is a linear mapping, then the map $d\mu:\mathfrak{J}\times \mathfrak{J} \to \mathfrak{J}$ given by
$$d\mu(a,b) = \mu(a)b + a\mu(b) - \mu(ab)$$
lies in $\mathbb{V}$. The quotient space of $\mathbb{V}$ by the space consisting of the maps $d\mu$ is called the second cohomology group of $\mathfrak{J}$ with coefficients in itself and is denoted by $H^2(\mathfrak{J},\mathfrak{J}).$

\begin{Th}
\label{marg_rig} The second cohomology group of the algebra $\mathfrak{J}_n$ is zero. Hence, $\mathfrak{J}_n$ is rigid.
\end{Th}
\begin{Proof} Let $h$ be any bilinear map from the space $\mathbb{V}.$ Let $a = b = e$ in the relation (\ref{cohom_2}). Reducing the identical terms that appear on both sides of the equation, we have
$$(h(e,e)e)e = h(e,e)e,$$
hence $h(e,e)e \in (\mathfrak{J}_n)_2.$ Using the multiplication table of $\mathfrak{J}_n$ it is now easy to see that 
\begin{equation}
\label{h(e,e)}
h(e,e) \in (\mathfrak{J}_n)_2.
\end{equation}
Let $a \in (\mathfrak{J}_n)_1, b = e$ in the relation (\ref{cohom_2}):
$$(h(a,a)e)a = \frac{1}{2}h(a,a)a$$
Representing the element $h(a,a)$ as the sum of its Peirce components, we have
$$((h(a,a)_1 + h(a,a)_2)e)a = \frac{1}{2}(h(a,a)_1 + h(a,a)_2)a.$$
Simplifying, we have $h(a,a)_2a = 0,$
which implies that $h(a,a)_2 = 0,$ since the Peirce component $(\mathfrak{J}_n)_2$ is one-dimensional. Linearizing the last relation (i.e., putting in it $a + b, b \in (\mathfrak{J}_n)_1$ instead of $a$) and using the relation (\ref{cohom_1}) we have
\begin{equation}
\label{h(a,b)1}
h(a,b) \in (\mathfrak{J}_n)_1, \text{ where } a, b \in (\mathfrak{J}_n)_1.
\end{equation}
Now let us linearize the relation (\ref{cohom_2}), that is, substitute the element $a$ in it with the sum $a+c.$ The resulting relation (which, for our convenience, we will write as $h = 0$) has 36 terms, but we can simplify it consideratebly. Write the expression $h$ as the sum of its ''homogeneous components'': $h = \sum_{i_a,i_b,i_c} h_{(i_a,i_b,i_c)},$ where $h_{(i_a,i_b,i_c)}$ is the sum of terms in $h$ in which the element $a$ occurs $i_a$ times, $b$ occurs $i_b$ times and $c$ occurs $i_c$ times (in other words, if $a, b$ and $c$ are considered as variables, then $h_{(i_a,i_b,i_c)}$ is the homogeneous component of $h$ of the multidegree $(i_a,i_b,i_c)$). Then it is easy to see that
$$h = h_{(3,1,0)}(a,b,c) + h_{(2,1,1)}(a,b,c) + h_{(1,1,2)}(a,b,c) + h_{(0,1,3)}(a,b,c) = 0.$$
The terms $h_{(3,1,0)}$ and $h_{(0,1,3)}$ are in fact zero, since they are exactly the left side of the relation (\ref{cohom_2}) (the term $h_{(0,1,3)}$ with the element $c$ instead of $a$). Hence we are left with the relation
\begin{equation}
\label{h_sum}
h_{(2,1,1)}(a,b,c) + h_{(1,1,2)}(a,b,c) = 0.
\end{equation}
Now let $0, 1  \neq \lambda   \in \mathbb{C}.$ Substituting in the relation (\ref{h_sum}) $\lambda a$ instead of $a,$ we have
$$\lambda^2h_{(2,1,1)}(a,b,c) + \lambda h_{(1,1,2)}(a,b,c) = 0.$$
Dividing the last relation by $\lambda$ and subtracting the relation (\ref{h_sum}), we have $h_{(2,1,1)} = 0.$ This relation has 12 terms and can be written explicitely:
\begin{equation}
\label{cohom_3}
2(h(a,c)b)a + (h(a,a)b)c + 2h(ac,b)a + h(a^2,b)c + 2h((ac)b,a) + h(a^2b,a) =
\end{equation}
$$2(ac)h(b,a) + a^2h(b,c) + 2h(a,c) + h(a,a)(bc) + 2h(ac,ba) + h(a^2,bc).$$
Let $a = b = e, c \in \mathfrak{J}_1$ in the relation (\ref{cohom_3}). Using the relation (\ref{h(e,e)}) and reducing the identical terms that appear on both sides of the equation, we have
$$2(h(e,c)e)e + \frac{1}{2}h(e,e)c = 2h(e,c)e.$$
Representing the element $h(e,c)$ as the sum of its Peirce components, we have
$$2h(e,c)_2 + \frac{1}{2}h(e,c)_1 + \frac{1}{2}h(e,e)c = 2h(e,c)_2 + h(e,c)_1.$$
Simplifying, we obtain
\begin{equation}
\label{h(e,c)}
h(e,e)c = h(e,c)_1, \text{ where } c \in (\mathfrak{J}_n) _1.
\end{equation}
Let $a \in (\mathfrak{J}_n)_1, b, c = e$ in the relation (\ref{cohom_3}). Using the relation (\ref{h(a,b)1}) and reducing the identical terms that appear on both sides of the equation, we have
$$2(h(a,e)e)a = h(a,e)a + \frac{1}{4}h(a,a).$$
Representing the element $h(a,e)$ as the sum of its Peirce components, we have
$$h(a,a) = 4h(a,e)_2a.$$
Linearizing this relation, we have
\begin{equation}
\label{h(a,b)}
h(a,b) = 2(h(a,e)_2b + h(b,e)_2a), \text{ where } a, b \in (\mathfrak{J}_k)_1.
\end{equation}
The relations (\ref{h(e,e)}), (\ref{h(e,c)}) and (\ref{h(a,b)}) are satisfied by any element $h$ of the space $\mathbb{V}$. 
Now we can show the statement of the proposition. To do this, we construct a map $\mu: \mathfrak{J} \to \mathfrak{J}$ such that $h = d\mu.$ Firstly, take a linear map $\mu_1$ defined by $\mu_1(a) = h(a,e).$ Then by (\ref{h(e,e)})
$$d\mu_1(e,e) = h(e,e)e + eh(e,e) - h(e,e) = h(e,e).$$
Thus for the map $h' = h - d\mu_1$ we have $h'(e,e) = 0$, therefore by (\ref{h(e,c)}) we have $h'(e,c)_1 = 0.$ Now let $\mu_2$ be a linear map on $\mathfrak{J}$ defined by $\mu_2(a) = 2h'(a,e)$. Then for $a, b \in \mathfrak{J}_1$ we have
$$d\mu_2(e,e) = 0 = h'(e,e),$$
$$d\mu_2(a,e) = 2h'(a,e)e + 2h'(e,e)a -h'(a,e) = \text{ (since $h'(a,e)_1 = 0$) } = h'(a,e).$$
$$d\mu_2(a,b) = 2h(a,e)b + 2h(b,e)a - 2h(ab,e) = \text{ (since $h'(a,e)_1 = h'(b,e)_1 = 0$) }$$
$$= 2h'(a,e)_2b + 2h'(b,e)_2a = \text{ (by (\ref{h(a,b)})) } = h'(a,b).$$
Hence $h' = d\mu_2$ and $h = d\mu_1 + d\mu_2.$ Therefore, $h = d\mu,$ where $\mu$ is the map given by
$$\mu(a) = 3h(e,a) -2(h(e,e)a + h(e,a)e + h(ea,e)),$$
and $H^2(\mathfrak{J}_n,\mathfrak{J}_n) = 0.$ 
\end{Proof}

\medskip
Now we see that the algebras $\mathfrak{J}_n$ stand out in a certain way. They are rigid, but can only degenerate to  zero algebra, which means that the closure of the ${\rm GL}_n$-orbit of the point corresponding to the algebra $\mathfrak{J}_n$ has a one-point intersection with the intersection of all other irreducible components of $\mathfrak{Jord}_n.$ Therefore, there exists no non-trivial Jordan algebra of dimension $n \geq 2$  which lies in the intersection of all irreducible components in the variety of $n$-dimensional Jordan algebras. Note that in the variety of Lie algebras such algebra exists \cite{gorb12}.

\subsection{Associative algebras}
The variety of associative algebras is defined by the following identity:
$$(xy)z=x(yz).$$
In the variety of $n$-dimensional associative algebras $\mathfrak{Ass}_n$
there are two marginal component defined by the following algebras:
$$ \nu^l_n = \langle e, n_1, \ldots, n_{n-1}\rangle \ : \ e^2=e, \ en_i=n_i \ (i=1, \ldots, n-1),   $$
$$ \nu^r_n = \langle e, n_1, \ldots, n_{n-1}\rangle \ : \ e^2=e, \ n_ie=n_i  \ (i=1, \ldots, n-1).   $$
As follows from below these algebras are rigid in a more bigger variety.

\subsection{Alternative algebras}
The variety of alternative algebras is a generalization of associative algebras and defined by the following identities:
$$x^2y=x(xy), \ xy^2=(xy)y.$$
In the variety of $n$-dimensional alternative algebras $\mathfrak{Alt}_n$
there are two marginal component defined by  $\nu^l_n$ and $\nu^r_n.$

\subsection{Left alternative algebras}
The variety of left alternative algebras is a generalization of associative and alternative algebras and defined by the following identity:
$$x^2y=x(xy).$$
In the variety of $n$-dimensional left alternative algebras $\mathfrak{LAlt}_n$
there are two marginal component defined by  $\nu^l_n$ and $\nu^r_n.$

Let us recall the definition of the second cohomology group of a left alternative algebra: 
Let $\mathfrak{L}$ be a left alternative  algebra. 
Let $\mathbb{V}$ be the space of the bilinear maps $h: \mathfrak{L}\times\mathfrak{L} \to \mathfrak{L}$ such that (see \cite{Jac})
\begin{equation}
\label{cohom_lefta}
h(a,a)b + h(a^2,b) = ah(a,b) + h(a,ab).
\end{equation}
For example, if $\mu: \mathfrak{L} \to \mathfrak{L}$ is a linear mapping, 
then the map $d\mu:\mathfrak{L}\times \mathfrak{L} \to \mathfrak{L}$ given by
$$d\mu(a,b) = \mu(a)b + a\mu(b) - \mu(ab)$$ lies in $\mathbb{V}$. 
The quotient space of $\mathbb{V}$ by the space consisting of the maps $d\mu$ is called the second cohomology group of $\mathfrak{L}$ with coefficients in itself and is denoted by $H^2(\mathfrak{L},\mathfrak{L}).$

Now, we will consider the algebra $\nu^l_n$ as a left alternative algebra and will prove that $\nu^l_n$ is rigid in the variety of 
$n$-dimensional left alternative algebras.
\begin{Th}
\label{marg_lalt} 
The second cohomology group of the left alternative algebra  $\nu^l_n$  is zero. Hence, $\nu^l_n$ is rigid.
\end{Th}
\begin{Proof} Let $h$ be any bilinear map from the space $\mathbb{V}.$ Let $a = b = e$ in the relation (\ref{cohom_lefta}). Reducing the identical terms that appear on both sides of the equation, we have
$h(e,e)e=h(e,e)$ and $h(e,e)\in (\nu^l_n)_{2}.$
Using the multiplication table of the algebra it is now easy to see that
for $a=e, b \in  (\nu^l_n)_{1}$ we have $h(e,e)b=h(e,b).$
It is follows that  $h(e,b)_1=0.$
On the other hand, for $a \in  (\nu^l_n)_{1}, b=e$ we can see that $h(a,a)e=0,$ and  $h(a,a)_1=0.$

From (\ref{cohom_lefta}), by linearization, we have:
\begin{equation}
\label{cohom_lefta_lin}
h(a,c)b+h(c,a)b+h(ac,b)+h(ca,b)=ah(c,b)+ch(a,b)+h(a,cb)+h(c,ab).\end{equation}

By some similar calculations, from (\ref{cohom_lefta_lin}),
for 
$$a=b=e, c\in  (\nu^l_n)_{1}; \ a \in  (\nu^l_n)_{1}, b=c=e; \ a \in  (\nu^l_n)_{1}, b=e, c \in  (\nu^l_n)_{1};$$ 
$$a \in  (\nu^l_n)_{1}, b \in  (\nu^l_n)_{1}, c=e; \ a,b,c \in  (\nu^l_n)_{1}$$ 
we have the following relations:
$$h(e,e)_1=0, \ h(e,a)_1=0, \ h(e,a)=h(e,e)a, \ h(a,e)_1=0, \ h(a,b)_1=0, \ h(a,b)=h(a,c)b.$$

Now we can show the statement of the proposition. 
To do this, we construct a map $\mu: \nu^l_n  \to \nu^l_n$ such that $h = d\mu.$ 

We take $\mu (x)=h(x,e).$ 
Obviously, we have

$$d\mu(e,e)=\mu(e)e=h(e,e), \ d\mu(a,e)=\mu(a)e=h(a,e)e=h(a,e),$$
$$d\mu(e,a)=\mu(e)a=h(e,e)a=h(e,a), \ d\mu(a,b)=\mu(a)b=h(a,e)b=h(a,b).$$

Now, $h=d\mu,$ hence $H^2(\nu^l_n,\nu^l_n)=0.$
\end{Proof}

By some similar calculations, we have the following 

\begin{Th}
\label{marg_lalt} 
The second cohomology group of the left alternative algebra  $\nu^r_n$  is zero. Hence, $\nu^r_n$ is rigid.
\end{Th}

\subsection{Non commutative Jordan algebras}
The variety of non commutative Jordan algebras is a generalization of associative, alternative and Jordan algebras and defined by the following identities: 
$$(xy)x=x(yx), \ (x^2y)x=x^2(yx).$$

\begin{Th}
\label{marg_flex} 
The non commutative Jordan algebras  $\nu^l_n, \nu^r_n$ and $\mathfrak{J}_n$ are not rigid. 
\end{Th}
\begin{Proof} 
Let us consider the family $\nu_n(\alpha)$ of algebras defined by
$$ \nu_n(\alpha) = \langle e, n_1, \ldots, n_{n-1}\rangle \ : \ e^2=e, \ en_i=\alpha n_i, \ n_ie=(1-\alpha) n_i  \ (i=1, \ldots, n-1; \alpha \in \mathbb{C}).   $$
For every $\alpha \in \mathbb{C}$ we have that $\nu_n(\alpha)$ is a non commutative Jordan algebra (it is easy to see, that $\nu_n(\alpha)^{(+)}=\mathfrak{J}_n$).
Now, the algebras $\nu^l_n, \nu^r_n$ and $\mathfrak{J}_n$ are in the orbit closure of the family $\nu_n(\alpha).$
\end{Proof}

\subsection{Leibniz algebras}
The variety of Leibniz algebras is a generalization of Lie algebras  and defined by the following identity:
$$(xy)z=(xz)y+x(yz).$$
We have a similar result in the variety of $n$-dimensional Leibniz algebras $\mathfrak{Leib}_n$ \cite{Ancochea,khom}.
A marginal irreducible component of $\mathfrak{Leib}_n$ consists of $n$-dimensional algebras with $(n-1)$-dimension left annihilator.
Obviously, every Lie algebra with  $(n-1)$-dimension left annihilator has zero multiplication.

\subsection{Anticommutative algebras}
There are many well known varieties of anticommutative algebras including Lie algebras:
\begin{enumerate}
\item[$\bullet$] Malcev algebras  \cite{kpv};
\item[$\bullet$] Binary Lie algebras \cite{kpv};
\item[$\bullet$] Sagle algebras \cite{fil};
\item[$\bullet$] Anti-commutative algebras with skew-symmetric identities \cite{dzhuma}
\end{enumerate}

and some others.

\begin{Th}
Let $\Omega$-algebras be a class of algebras including all Lie algebras and satisfying 
 the family of multilinear polinomial identities $\Omega$,
then the intersection of all irreducible components of the variety  $\mathfrak{AVar}_n$ of $n$-dimensional anticommutative $\Omega$-algebras
has an algebra with non zero multiplication.

\end{Th}

\begin{Proof}
In the variety of $2$-dimensional algebras there is only one anticommutative algebras with nonzero multiplication. 
It is the solvable (non-nilpotent) Lie algebra with multiplication $e_1e_2=-e_2e_1=e_2.$
In \cite[Corollary 5.9]{kv172} was proved that every anticommutative algebra with level two is a Lie algebra.
Now, from the list of algebras in \cite[Corollary 5.9]{kv172}, it is easy to see, that
every anticommutative algebra of level two degenerates to the nilpotent anticommutative algebra $\mathfrak{n}_3$  with the multiplication $e_1e_2=-e_2e_1=e_3.$
It is follows that $\mathfrak{n}_3$ is in the intersection of all irreducible components of $\mathfrak{AVar}_n.$
Now we have the statement of the Theorem. 
\end{Proof}


\end{document}